\documentclass[smallheadings,11pt]{scrartcl}
\usepackage{amsmath}
\usepackage{amssymb}
\usepackage{avant}
\usepackage[slantedGreek]{mathpazo}
\usepackage[only,llbracket,rrbracket]{stmaryrd}

\usepackage{graphicx}
\usepackage{color}
\usepackage[sort&compress,numbers]{natbib}
\usepackage{fancyhdr}
\usepackage{bm}
\usepackage[T1]{fontenc}
\newcommand{\disc}[1]{\left\llbracket {#1} \right\rrbracket}
\newcommand{\ave}[1]{\left< {#1} \right>}

\newcommand{\brac}[1]{\left( {#1} \right)}
\newcommand{\bracc}[1]{\left\{ {#1} \right\}}

\newcommand{\bs}[1]{\bm{#1}}
\newcommand{\vect}{\bs}
\newcommand{\pd}{\partial}
\newcommand{\dif}{ \ d}

\newcommand{\ip}{\mathbf{:}}
\newcommand{\grad}{{\nabla}^{\rm s}}
\newcommand{\gradn}{{\nabla}}

\newcommand{\eps}{\epsilon}
\newcommand{\strain}{\bs{\eps}}
\newcommand{\stress}{\bs{\sigma}}
\newcommand{\trans}[1]{ {#1}^{\text{T}} }

\usepackage{time}

\newcommand{\epsbc}{\eps_{\rm bc}}

\renewcommand{\Bar}{\overline}
\renewcommand{\Tilde}{\widetilde}

\usepackage{sectsty}
\sectionfont{\large}
\subsectionfont{\normalsize}
\usepackage{setspace}
\title{	  \Large  A discontinuous Galerkin formulation for a strain 
					gradient-dependent damage model}
\author{\normalsize Garth N. Wells$^{1}$\footnote{Corresponding author, email:
g.n.wells@citg.tudelft.nl, fax: +31 15 278 6383.} \ \ \  
			Krishna Garikipati$^{2}$  \ \ \  
			Luisa Molari$^{3}$}

\date{\small 
		$^{1}$Faculty of Civil Engineering and Geosciences, 
		Delft University of Technology \\ Stevinweg 1, 2628 CN Delft, The
		Netherlands \\[1ex]
		$^{2}$Department of Mechanical Engineering,
		University of Michigan \\ Ann Arbor, Michigan 48109, USA \\[1ex]
		$^{3}$DISTART, Universit\`{a} di Bologna \\
			Viale Risorgimento 2, 40136 Bologna, Italy \\[2ex]
}
\begin{document}

\maketitle
\onehalfspacing

\begin{abstract}
	\noindent  The numerical solution of strain gradient-dependent continuum 
	problems
	has been dogged by continuity demands on the basis functions. For most
	commonly  accepted models, solutions using the finite element method 
	demand $C^{1}$ continuity of the shape functions. 
	Here, recent development in
	discontinuous Galerkin methods are explored and exploited for the
	solution of a prototype nonlinear strain gradient dependent continuum model.
	A formulation is developed that allows the rigorous solution of a strain 
	gradient
	damage model using standard $C^{0}$ shape functions.
	The
	formulation is tested in one-dimension for the simplest possible
	finite element formulation: piecewise linear displacement and constant 
	(on elements) internal variable. Numerical results are shown to compare 
	excellently with a benchmark solution. The results are remarkable given the
	simplicity of the proposed formulation.

\end{abstract}
\newpage
\doublespacing
\section{Introduction}
	Strain gradient dependent continuum models have been developed for a wide
	range of problems. Strain gradient effects are included in continuum models
	to reproduce experimentally observed phenomena which cannot be 
	captured with classical
	models~\citep{aifantis:1984,larsy:1988,muhlhaus:1991,peerlings:1996,fleck:1997,
	gao:1999,fleck:2001}.  
	The range of application is broad, from
	large geological problems to polycrystals. 
	Typical phenomena which can be captured with strain gradient models include
	strain localisation in the presence of softening and size effects.
	
	The development of strain gradient models has
	been hindered by the lack of a suitable numerical framework
	 for their robust solution on arbitrary domains.
	The introduction of strain gradients into continuum models poses
	significant challenges in solving the ensuing equations. The finite element
	method, the dominant numerical method in solid mechanics, is
	ideally suited to the solution of second-order partial differential
	equations, such as classical elasticity. The solution of
	gradient-dependent continuum problems usually demands at least 
	$C^{1}$ continuous
	basis functions, which are difficult to construct in spatial 
	dimensions higher than one. 
	Previous attempts to solve such problems with $C^{1}$ shape
	functions or \emph{ad-hoc} measures have proven 
	difficult~\citep{larsy:1988,borst:1992,borst:1996}.
	More seriously, in numerous publications, basic 
	continuity requirements are completely ignored.  
	To avoid these difficulties, \citet{askes:2000b}
	applied the element-free Galerkin method, which can provide a high degree of
	continuity, for the solution of strain gradient dependent damage models.
	However, the element-free Galerkin method entails other difficulties, 
	lacks the penetration  in the solid mechanics community 
	of the finite element method, and is generally less efficient. 
	As a result of these difficulties, strain gradient dependent models 
	are not widely applied, and many formulations are largely untested.
	The difficulties presented by continuity requirements has even 
	lead to reformulations of strain gradient models that are driven by 
	algorithmic convenience~\citep{peerlings:1996,engelen:2003}.

	In this work, a fresh perspective is taken on the solution of
	strain gradient dependent continuum problems in light of recent developments
	in discontinuous and continuous/discontinuous Galerkin methods for elliptic
	problems~\citep{oden:1998dg,arnold:2002,engel:2002}. A summary of recent
	developments can be found in~\citet{arnold:2002}.
	In the derivation of the Galerkin
	problem, potential discontinuities in the basis functions 
	across internal surfaces  are taken into account, resulting in a
	generalisation of the conventional Galerkin method.

	To begin, a strain gradient-dependent damage model, which is used as a 
	prototype example, is introduced. It is cast as a continuous 
	Galerkin problem in a finite element framework 
	and the difficulties with the conventional finite element method 
	are highlighted. The Galerkin problem is
	then generalised to allow for discontinuities in the appropriate fields.
	The formulation is tested for the simplest possible finite element 
	in one dimension. A series of test cases are computed and the results
	are compared to a benchmark solution.

\section{Gradient-enhanced damage model: Preliminaries}	
	Consider a body $\Omega$ in $\mathbb{R}^{n}$, with boundary $\Gamma = \pd
	\Omega$.
	The strong form of the equilibrium equation for the body $\Omega$, in the
	absence of body forces, and 
	associated standard boundary conditions, is:
	\begin{align}
		\gradn \cdot \stress 	&= \vect{0}  &{\rm in} \ \Omega 
					\label{eqn:strong_equil} \\
		\stress \cdot \vect{n} 	&= \vect{h}  &{\rm on} \ \Gamma_{h} \\
		\vect{u} 				&= \vect{g}  &{\rm on} \ \Gamma_{g} 
	\end{align}		
	where $\gradn$ is the gradient operator,
	$\stress$ is the stress tensor, $\vect{h}$ is the prescribed traction
	on $\Gamma_{h}$ and $\vect{g}$ is the prescribed displacement on the
	boundary $\Gamma_{g}$ ($\Gamma_{g} \cup \Gamma_{h} = \Gamma$, 
	$\Gamma_{g} \cap \Gamma_{h} = \emptyset$). The outward normal to $\Gamma$ is
	denoted~$\vect{n}$.

	For an isotropic elasticity-based damage model, 
	the stress at a material point is given by:
	\begin{equation}
		\stress = \brac{1-\omega} \mathcal{C} \ip \grad \vect{u} 
	\label{eqn:constit}
	\end{equation}	
	where $\mathcal{C}$ is the usual linear-elastic constitutive tensor 
	and the damage variable ($0 \le \omega \le 1$) is a function of a
	scalar history parameter~$\kappa$,
	\begin{equation}
		\omega = \omega\brac{\kappa}. 
	\label{eqn:damage}
	\end{equation}	
	The history parameter $\kappa$ is related to a gradient-dependent 
	`equivalent strain', $\Bar{\eps}$.
	A common choice for $\Bar{\eps}$ is: 
	\begin{equation}
		\Bar{\eps} = \eps_{\rm eq} + c^{2} \Delta \eps_{\rm eq} 
	\label{eqn:grad_strain}
	\end{equation}	
	where $\eps_{\rm eq}$ is an invariant of the local strain tensor 
	$\strain =\grad \vect{u}$, 
	$c$ is a length scale which reflects the strength of
	strain gradient effects and $\Delta$ is the Laplacian operator.  
	This formulation is often named `explicit gradient 
	damage'~\citep{peerlings:1996}.
	The chosen invariant for the local equivalent
	strain reflects the processes that drive damage growth in a given 
	material. In one  dimension,
	the obvious choice is that the equivalent strain is equal to the strain. 

	The history parameter $\kappa$ is equal to the largest positive value of 
	$\Bar{\eps}$ reached at a material point. Defining a loading function $f$,
	\begin{equation}
		f = \Bar{\eps} - \kappa
	\end{equation}	
	the evolution of $\kappa$ obeys the Kuhn-Tucker conditions,
	\begin{align}
		\Dot{\kappa} \ge 0, && f \le 0,  && \Dot{\kappa}f = 0.
	\end{align}	
	 A commonly adopted dependency is:
	\begin{equation}
		\omega = 
		\begin{cases}
			0			& {\rm if} \ \ \kappa \le \kappa_{0} \\
			1- \dfrac{\kappa_{0} \brac{\kappa_{c} - \kappa}}{\kappa
				\brac{\kappa_{c}-\kappa_{0}}} 
						& {\rm if} \ \ \kappa_{0} < \kappa < \kappa_{c} \\
			1			& {\rm if} \ \ \kappa \ge  \kappa_{c}
		\end{cases}		
	\label{eqn:damage_evolution}
	\end{equation}	
	where $\kappa_{0}$ is the value of the history parameter 
	at which damage begins to develop  and
	$\kappa_{c}$ is the value at which $\omega=1$.
	The evolution of $\omega$ in equation~\eqref{eqn:damage_evolution} 
	yields a linear softening response for a uniaxial
	test in the absence of strain gradient effects. To make the dependency of
	$\omega$  on $\Bar{\eps}$ clear, the expressions $\omega\brac{\kappa}$ and 
	$\omega\brac{\Bar{\eps}}$ will be used interchangeably.
	
	Insertion of the constitutive model (see equations~\eqref{eqn:constit},
	\eqref{eqn:damage} and~\eqref{eqn:grad_strain}) into the equilibrium 
	equation~\eqref{eqn:strong_equil} leads to 
	a non-linear fourth-order partial differential equation. This requires
	the prescription of boundary conditions on gradients of the displacement
	field higher than one. The physical implications of these boundary
	conditions are unclear and are the subject of debate. At this stage,
	the boundary condition
	\begin{align}
		c^{2} \gradn \eps_{\rm eq} \cdot \vect{n} 
					= \epsbc  && {\rm on} \ \Gamma
	\label{eqn:eps_bs}
	\end{align}	
	is considered. A common choice is $\epsbc =0$,
	which is adopted for all examples in Section~\ref{sec:examples}.

	This elasticity-based damage model is convenient for preliminary 
	developments as $\eps_{\rm eq}$
	is calculated explicitly from the gradient of the 
	displacement field, which is  in
	contrast to the equivalent plastic strain in an elastoplastic model. 
	However, equation~\eqref{eqn:grad_strain} is identical in form to the
	equation for the gradient-dependent equivalent plastic strain that is 
	adopted in many strain gradient dependent plasticity 
	models~\citep{aifantis:1984,muhlhaus:1991,fleck:2001}. 
	This model therefore provides a
	canonical formulation which can be extended to a broader class of models.

\section{Galerkin formulation}	
	In developing a weak formulation for eventual finite element solution, 
	the equilibrium equation~\eqref{eqn:strong_equil}  and the
	equation for $\Bar{\eps}$~\eqref{eqn:grad_strain}  are considered 
	separately. The non-linear fourth-order equation resulting from insertion 
	of  the constitutive equations into the
	equilibrium equation could potentially be cast in a weak from. 
	The formulation would inevitably be specific to the
	chosen dependency of damage on $\Bar{\eps}$, a dependency which is
	potentially highly complex. Hence, for simplicity and generality, it is
	convenient to treat the two equations separately.

	The body $\Omega$ is partitioned into $n_{el}$ non-overlapping elements
	$\Omega_{e}$ such that 
	\begin{equation}
		\Bar{\Omega} = \bigcup_{e=1}^{n_{el}} \Bar{\Omega}_{e}.
	\end{equation}	
	where $\Bar{\Omega}_{e}$ is a closed set (i.e., it includes the boundary of
	the element).
	The elements $\Omega_{e}$ (which are open sets) 
	satisfy the standard requirements for a finite  element partition. A
	domain $\Tilde{\Omega}$ is also defined
	\begin{equation}
		\Tilde{\Omega} = \bigcup_{e=1}^{n_{el}} \Omega_{e}
	\end{equation}	
	where $\Tilde{\Omega}$ does not include element boundaries. It is also
	useful to define the `interior' boundary $\Tilde{\Gamma}$,
	\begin{equation}
		\Tilde{\Gamma} = \bigcup_{i=1}^{n_{b}} \Gamma_{i}
	\end{equation}	
	where $\Gamma_{i}$ is the $i$th interior element boundary and $n_{b}$ is
	the number of internal inter-element boundaries.
	
	Consider now the function spaces $\mathcal{S}^{h}$, $\mathcal{V}^{h}$ and 
	$\mathcal{W}^{h}$,
	\begin{align}
		\mathcal{S}^{h} &=\bracc{u^{h}_{i} \in H_{0}^{1}\brac{\Omega} \ \left| \
				u_{i}^{h}|_{\Omega_{e}} \in P_{k_{1}}\brac{\Omega_{e}} \forall e,  
						\ u_{i} = g_{i} \ {\rm on}
						\ \Gamma_{g} \right. }	\label{eqn:disp_trial}	\\ 
		\mathcal{V}^{h} &=\bracc{w^{h}_{i} \in H_{0}^{1}\brac{\Omega} \ \left| \
				w_{i}^{h}|_{\Omega_{e}} \in P_{k_{1}}\brac{\Omega_{e}} \forall e,  
						\ w_{i} = 0 \ {\rm on}
						\ \Gamma_{g} \right. }	\label{eqn:disp_test}		\\
		\mathcal{W}^{h} &=\bracc{q^{h} \in L_{2}\brac{\Omega} \ \left| \
				q^{h}|_{\Omega_{e}} \in P_{k_{2}}\brac{\Omega_{e}} \forall e
				\right. }		\label{eqn:space_q}	
	\end{align}	 
	where $P_{k}$ represents the space of polynomial finite element shape 
	functions (of polynomial order $k$).
	The spaces $\mathcal{S}^{h}$ and $\mathcal{V}^{h}$ represent 
	usual $C^{0}$ continuous finite element shape functions. The space
	$\mathcal{W}^{h}$ can contain discontinuous functions.

\subsection{Standard Galerkin weak form}
	The standard, continuous Galerkin problem for the 
	equilibrium equation~\eqref{eqn:strong_equil} is of the form: Find
	$\vect{u}^{h} \in \mathcal{S}^{h}$ such that
	\begin{align}
		\int_{\Omega} \gradn \vect{w}^{h} \ip 
			\brac{1 - \omega\brac{\Bar{\eps}^{h}}} \mathcal{C} \ip 
					\grad \vect{u}^{h} \dif \Omega
		 - \int_{\Gamma_{h}} \vect{w}^{h} \cdot \vect{h} \dif \Gamma = 0
			&& \forall \vect{w}^{h} \in \mathcal{V}^{h} \label{eqn:equil_cont}
	\end{align}	
	where it was already assumed that $\vect{u}^{h}$ is $C^{0}$ 
	continuous (see equation~\eqref{eqn:disp_trial}).
	Note that the damage is a function of $\Bar{\eps}$, which is in turn 
	a function of
	displacement gradients, making the equation non-linear.
	It is presumed at this point that $\Bar{\eps}^{h}$ is square-integrable
	over $\Omega$ ($\Bar{\eps}^{h} \in L_{2}\brac{\Omega}$). 

	A second Galerkin problem is constructed to solve for $\Bar{\eps}$ 
	(equation~\eqref{eqn:grad_strain}). It consists of: Find
	$\Bar{\eps} \in \mathcal{W}^{h}$ such that  
	\begin{equation}
		\int_{\Omega} q^{h} \Bar{\eps}^{h} \dif \Omega
			-\int_{\Omega} q^{h} \eps_{\rm eq}^{h} \dif \Omega
			+\int_{\Omega} \gradn q^{h}\cdot c^{2} \gradn \eps_{\rm eq}^{h} 
				\dif \Omega 
			-\int_{\Gamma} q^{h} \epsbc
				\dif \Gamma 
			= 0 \ \ \ \  \forall q^{h} \in \mathcal{W}^{h}
	\label{eqn:bar_eps_weak}
	\end{equation}	
	where it is assumed that $\eps_{\rm eq}^{h}$ is known. Recall that
	discontinuities in $q^{h}$ and $\Bar{\eps}^{h}$ are permitted.

	Two difficulties exist in the preceding Galerkin formulation. The first 
	is  that the
	weight function $q^{h}$ can be discontinuous (cf.
	equation~\eqref{eqn:space_q}), meaning that $\gradn q^{h}$ is not
	necessarily square-integrable on $\Omega$. 
	This problem can
	be circumvented easily by requiring $C^{0}$ continuity of the  
	functions in~$\mathcal{W}^{h}$.	The second problem, which is 
	less easily solved, is 
	that $\eps_{\rm eq}^{h}$ is computed from $\grad \vect{u}^{h}$. 
	Therefore, calculating $\gradn \eps_{\rm eq}^{h}$ everywhere in $\Omega$
	requires that the displacement field $\vect{u}^{h}$ be $C^{1}$ continuous 
	if  singularities are to be avoided. However, since $\vect{u}^{h} \in
	H^{1}_{0}\brac{\Omega}$ (see equation~\eqref{eqn:disp_trial}), it is not
	necessarily $C^{1}$ continuous.
 
 	To proceed with this formulation in a conventional manner, 
	two possibilities present themselves.  
	The first is to solve
	equations~\eqref{eqn:equil_cont} and~\eqref{eqn:bar_eps_weak} using $C^{0}$
	finite element shape functions to interpolate 
	$\Bar{\eps}^{h}$ and $q^{h}$, which is straightforward, and using $C^{1}$
	shape functions for $\vect{w}^{h}$ and $\vect{u}^{h}$. 
	The second approach
	is to interpolate $\eps_{\rm eq}$ using $C^{1}$ shape functions, from which
	the term $\Delta \eps_{\rm eq}$ can be evaluated everywhere in~$\Omega$.
	The second approach may appear more attractive than the first as it
	requires a $C^{1}$ interpolation of a scalar field rather than a vector
	field.
	Both approaches pose significant difficulties as $C^{1}$ shape functions are
	difficult to construct, lack generality and lead to extremely complex element
	formulations. $C^{1}$ functions are difficult to construct in two
	dimensions, and to the authors' knowledge, untried in three dimensions.

\subsection{Discontinuous Galerkin form}	 
	The approach advocated here avoids the need for $C^{1}$ continuity of the 
	displacement field by imposing the required degree of continuity in a weak
	sense. 
	Before proceeding with the formulation, it is necessary to define 
	jump and an averaging operations.
	The jump in a field $\vect{a}$ across a surface (which is associated with a
	body) is given by~\citep{arnold:2002}:
	\begin{equation}
		\disc{\vect{a}} = \vect{a}_{1} \cdot \vect{n}_{1} + 
		\vect{a}_{2} \cdot \vect{n}_{2} 
	\end{equation}		
	where the subscripts denote the side of the surface and $\vect{n}$ is the
	outward unit normal vector. This definition is 
	convenient as it avoids introducing `+' and `-' sides of a surface.
	This is particularly so for arbitrarily-oriented surfaces in two and
	three dimensions.	
	The average of a field
	$\vect{a}$ across a surface is given by:
	\begin{equation}
		\ave{\vect{a}} = \frac{\brac{\vect{a}_{1} + \vect{a}_{2}}}{2}.
	\end{equation}		

	Consider now equation~\eqref{eqn:grad_strain} for $\Bar{\eps}$, which can be
	cast in a weak form using integration by parts and 
	the divergence theorem on the boundary $\Gamma$ and on inter-element 
	boundaries,~$\Tilde{\Gamma}$. This yields:
	\begin{multline} \label{eqn:dg_1}
		\int_{\Omega} q^{h} \Bar{\eps}^{h} \dif \Omega
			-\int_{\Omega} q^{h} \eps_{\rm eq}^{h} \dif \Omega
			+\int_{\Tilde{\Omega}} \gradn q^{h} \cdot  c^{2}\gradn \eps_{\rm eq}^{h} 
					\dif \Omega 
			-\int_{\Gamma} q^{h} \epsbc \dif \Gamma \\
			-\int_{\Tilde{\Gamma}} \ave{q^{h}}  
						\cdot c^{2}\disc{ \gradn \eps_{\rm eq}^{h}} \dif \Gamma
			-\int_{\Tilde{\Gamma}} \disc{q^{h}} \cdot c^{2}\ave{\gradn \eps_{\rm
			eq}^{h}} \dif \Gamma
			= 0.
	\end{multline}	
	Note the distinction between $\Omega$ and~$\Tilde{\Omega}$ for the volume
	integrals.
	It is chosen that the following weak statements of continuity should hold:
	\begin{align}
		 \int_{\Tilde{\Gamma}} \ave{q^{h}} c^{2} \disc{\gradn \eps_{\rm eq}^{h}}
						\dif \Gamma &= 0  
						&& \forall q^{h} \in \mathcal{W}^{h} \\
		- \int_{\Tilde{\Gamma}} \ave{\gradn q^{h}} \cdot  c^{2} \disc{\eps_{\rm eq}^{h}}
						\dif \Gamma &= 0 
						&& \forall q^{h} \in \mathcal{W}^{h}. \label{eqn:dg_symm}
	\end{align}	
	Also, a `penalty-like' term is introduced:
	\begin{equation}
			 \int_{\Tilde{\Gamma}} \frac{c^{2}}{h_{e}}
			\disc{q^{h}}\cdot{\disc{\eps_{\rm eq}^{h}}} \dif \Gamma	=0	
	\end{equation}	
	where $h_{e}$ is a length scale which is required for dimensional 
	consistency. 
	Adding the additional equations to equation~\eqref{eqn:dg_1} 
	leads to the following Galerkin 
	problem: Find $\Bar{\eps}^{h}  \in \mathcal{W}^{h}$ such that		
	\begin{multline}
		\int_{\Omega} q^{h} \Bar{\eps}^{h} \dif \Omega
			-\int_{\Omega} q^{h} \eps_{\rm eq}^{h} \dif \Omega
			+\int_{\Tilde{\Omega}} \gradn q^{h} \cdot  c^{2}\gradn 
			\eps_{\rm eq}^{h} \dif \Omega 
			-\int_{\Gamma} q^{h} \epsbc \dif \Gamma \\
			-\int_{\Tilde{\Gamma}} \disc{q^{h}} \cdot 
					c^{2}\ave{\gradn \eps_{\rm eq}^{h}} \dif \Gamma 
			-\int_{\Tilde{\Gamma}} \ave{\gradn q^{h}} \cdot 
						c^{2}\disc{\eps_{\rm eq}^{h}} \dif \Gamma \\
			+ \int_{\Tilde{\Gamma}} \frac{c^{2}}{h_{e}}
			\disc{q^{h}}\cdot{\disc{\eps_{\rm eq}^{h}}} \dif \Gamma
			= 0 \ \ \ \ \forall q^{h} \in \mathcal{W}^{h}.
		\label{eqn:eps_weak}
	\end{multline}	
	Adding the term in equation~\eqref{eqn:dg_symm} to the problem
	provides a degree of `symmetry' with the term $\int_{\Tilde{\Gamma}} 
	\disc{q^{h}}\cdot c^{2}\ave{\gradn \eps_{\rm eq}^{h}} \dif \Gamma$. 
	The choice of $c^{2}/h_{e}$ may seem somewhat arbitrary
	considering that it appears as a penalty-like parameter. This choice will 
	be
	justified later through an analogy between the proposed method and
	a finite difference scheme. No gradients of $\eps_{\rm eq}^{h}$ or
	$q^{h}$ appear in terms integrated over $\Omega$ (which includes interior
	boundaries) in equation~\eqref{eqn:eps_weak}, hence the continuity
	requirements on the spaces $\mathcal{S}^{h}$ and $\mathcal{W}^{h}$ are
	sufficient.
	 
	Equation~\eqref{eqn:eps_weak} reassembles the
	`interior penalty' method for classical elasticity, which belongs to 
	the discontinuous Galerkin family of methods~\citep{arnold:2002}.
	Terms have been added to the weak form that for a
	conventional elasticity problem would lead to a symmetric formulation. 
	Symmetry is however not
	of relevance here as the functions $q^{h}$ and $\eps_{\rm eq}^{h}$ will 
	generally come from different function spaces. This formulation is 
	general for the
	case in which the space $\mathcal{W}^{h}$ contains discontinuous functions. 
	However, note 
	if all functions in the space $\mathcal{W}^{h}$ are 
	$C^{0}$ continuous, the formulation is still valid, with
	terms relating to the jump in $\eps_{\rm}^{h}$ remaining. The
	formulation would then resemble a continuous/discontinuous 
	Galerkin method~\citep{engel:2002}.

	The solution of the gradient enhanced damage problem requires
	the simultaneous solution of equations~\eqref{eqn:equil_cont} 
	and~\eqref{eqn:eps_weak}, which are coupled. 
	In summary, the problem is: Find $\vect{u}^{h} \in
	\mathcal{S}^{h}$ and $\Bar{\eps}^{h} \in \mathcal{W}^{h}$ such that
	\begin{align}
		\int_{\Omega} \gradn \vect{w}^{h} \ip 
			\brac{1 - \omega\brac{\Bar{\eps}^{h}}} \mathcal{C} \ip 
					\grad \vect{u}^{h} \dif \Omega
		 - \int_{\Gamma_{h}} \vect{w}^{h} \cdot \vect{h} \dif \Gamma = 0
			&& \forall \vect{w}^{h} \in \mathcal{V}^{h} \label{eqn:equil_cont_b}
	\end{align}	
	\begin{multline}
		\int_{\Omega} q^{h} \Bar{\eps}^{h} \dif \Omega
			-\int_{\Omega} q^{h} \eps_{\rm eq}^{h} \dif \Omega
			+\int_{\Tilde{\Omega}} \gradn q^{h} \cdot  c^{2}\gradn \eps_{\rm eq}^{h} 
					\dif \Omega 
			-\int_{\Gamma} q^{h} \epsbc \dif \Gamma \\
			-\int_{\Tilde{\Gamma}} \disc{q^{h}} \cdot 
					c^{2}\ave{\gradn \eps_{\rm eq}^{h}} \dif \Gamma 
			-\int_{\Tilde{\Gamma}} \ave{\gradn q^{h}} \cdot 
						c^{2}\disc{\eps_{\rm eq}^{h}} \dif \Gamma \\
			+ \int_{\Tilde{\Gamma}} \frac{c^{2}}{h_{e}}
			\disc{q^{h}}\cdot{\disc{\eps_{\rm eq}^{h}}} \dif \Gamma
			= 0 \ \ \ \ \forall q^{h} \in \mathcal{W}^{h}
		\label{eqn:eps_weak_b}
	\end{multline}	
	where the nonlinear equations are coupled through the dependency of 
	$\omega$ on $\Bar{\eps}^{h}$ and the dependency of $\Bar{\eps}^{h}$ 
	on~$\vect{u}^{h}$. Linearisation of these equations is straightforward, and
	is included in Appendix~\ref{append:linearisation}.

	In this work, the simplest possible finite element formulation is 
	considered. It is chosen to interpolate the displacement field with linear 
	piecewise continuous ($C^{0}$) functions and to use constant functions on 
	elements for $\Bar{\eps}$ ($k_{2} = 0$ in equation~\eqref{eqn:space_q}). 
	Also, the boundary condition $\gradn \eps_{\rm eq} \cdot
	\vect{n} =0$ on $\Gamma$ is applied.
	As a consequence, several terms disappear from 
	equation~\eqref{eqn:eps_weak_b}, leading to the problem: Find $\Bar{\eps}^{h}
	\in \mathcal{W}^{h}$ such that
	\begin{align}
		\int_{\Omega} q^{h} \Bar{\eps}^{h} \dif \Omega
			-\int_{\Omega} q^{h} \eps_{\rm eq}^{h} \dif \Omega
			+ \int_{\Tilde{\Gamma}} \frac{c^{2}}{h_{e}}
			\disc{q^{h} }\cdot{\disc{\eps_{\rm eq}^{h}}} \dif \Gamma
			= 0 &&\forall q^{h} \in \mathcal{W}^{h}.
		\label{eqn:eps_weak_fe}
	\end{align}	
	If $c=0$, $\Bar{\eps} = \eps_{\rm eq}$ at all points in $\Omega$,  and the model
	reduces to a local damage formulation (no gradient effects). 
	In one dimension, $h_{e}$ is
	taken as~$\ave{h_{e}}$. A higher-dimension generalisation would be the
	distance between the centroid of the neighbouring elements. 
		
	A physical interpretation of equation~\eqref{eqn:eps_weak_fe} is 
	simple. The stronger the spatial variation in the strain field , 
	the larger the
	jumps in the strain across element boundaries. 
	Equation~\eqref{eqn:eps_weak_fe} sets $\Bar{\eps}$ equal
	to the local equivalent strain, and subtracts a component which is
	proportional to the equivalent strain jump and the material parameter
	$c^{2}$, effectively decreasing $\Bar{\eps}$ (relative the $\eps_{\rm eq}$)
	in the presence of rapid spatial variation in the strain field, which is
	manifest in the form of jumps in the strain across element boundaries.
	
	In practice,
	this finite element formulation is very simple. An element has three nodes.
	Displacement degrees of freedom are located at the two end-nodes, and a
	degree of freedom for $\Bar{\eps}$ is located at the centre node of each
	element. The standard loop over all elements in a mesh is performed, and in
	addition all interior interfaces are looped over. Despite the node
	corresponding to $\Bar{\eps}^{h}$ being internal to an element, it 
	cannot be eliminated at the element level.
	The element stiffness matrices for this formulation are elaborated in
	Appendix~\ref{append:fe}.

\subsection{Consistency of the discontinuous formulation}
	Having added non-standard terms to the weak form, it is important to prove
	consistency of the method. Applying integration by parts to the 
	integral over $\Tilde{\Omega}$ in equation~\eqref{eqn:eps_weak}, 
	\begin{multline}
		\int_{\Tilde{\Omega}} \gradn q^{h} \cdot c^{2} \gradn \eps_{\rm ep}^{h} \dif \Omega 
		=
			-\int_{\Tilde{\Omega}} q^{h} c^{2}\Delta \eps_{\rm eq}^{h} \dif \Omega
			+\int_{\Gamma} q^{h} c^{2} \gradn \eps_{\rm eq} \cdot \vect{n} \dif \Gamma  \\
			+\int_{\Tilde{\Gamma}} \ave{q^{h}}
					 c^{2} \disc{\gradn \eps_{\rm eq}} \dif \Gamma
			+\int_{\Tilde{\Gamma}} \disc{q^{h}} \cdot
					 c^{2} \ave{\gradn \eps_{\rm eq}} \dif \Gamma.
	\end{multline}	
	Inserting this expression into the infinite-dimensional version of 
	equation~\eqref{eqn:eps_weak}, and employing standard variational arguments,
	the following  Euler-Lagrange equations can be identified:
	\begin{align}
		\Bar{\eps} - \eps_{\rm eq} - c^{2} \Delta \eps_{\rm eq} &= 0 
				&& {\rm in} \ \ \Tilde{\Omega} \label{eqn:EL_a}\\
		 c^{2} \disc{\eps_{\rm eq}} 							&= 0 
				&& {\rm on} \ \ \Tilde{\Gamma} 
		 				\label{eqn:EL_b} \\ 			
		 c^{2} \disc{\gradn \eps_{\rm eq}} 						&= 0
				&& {\rm on} \ \ \Tilde{\Gamma}  \label{eqn:EL_c} \\ 
		 c^{2} \gradn \eps_{\rm eq} \cdot  \vect{n} 			&= \epsbc
				&& {\rm on} \ \ \Gamma \label{eqn:EL_d} 
	\end{align}
	Equation~\eqref{eqn:EL_a} is the original problem over element interiors 
	(see equation~\eqref{eqn:strong_equil}). Equations~\eqref{eqn:EL_b}
	and~\eqref{eqn:EL_c} impose continuity of the corresponding fields across
	element boundaries and equation~\eqref{eqn:EL_d} imposes the natural
	boundary condition on $\gradn \eps_{\rm eq} \cdot \vect{n}$.
	The Galerkin form (equation~\eqref{eqn:eps_weak}) can therefore be seen as
	the weak imposition of these Euler-Lagrange equations.

\subsection{Finite difference analogy}
	In one-dimension for equally spaced nodal points, it can be shown that the
	proposed formulation ($C^{0}$ linear $u^{h}$ and piecewise constant
	$\Bar{\eps}^{h}$) is equivalent to a finite difference
	scheme for calculating $\eps_{{\rm eq},xx}$ (which is equal 
	to $u_{,xxx}$  for $\eps_{\rm eq} = u_{,x}$) 
	at the centre of each element.  

\begin{figure}
		\center\includegraphics{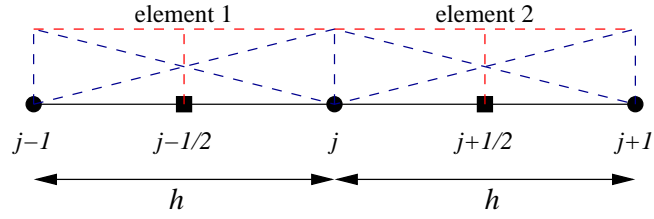}
	\label{fig:two_elements}	
	\caption{Two element configuration. Displacement degrees of freedom are
		located at the circular nodes, and $\Bar{\eps}$ degrees
		of freedom are located at the squares.} 
	\end{figure}

	Consider the two element configuration in figure~\ref{fig:two_elements}.
	The displacement degrees of freedom  are stored at the circular nodes and 
	are denoted~$a_{j}$. 
	From the form of the finite element shape functions, the jump in 
	the equivalent strain at element boundary $j$  is given by: 
	\begin{equation}
		\left. \frac{1}{h}\disc{\eps^{h}_{\rm eq}} \right|_{j} 
			= - \frac{a_{j-1} - 2a_{j} + a_{j+1}}{h^{2}}
			= - u^{\prime \prime}|_{j}
	\end{equation}
	which is equivalent to the second-order finite difference expression for the
	second derivative of the displacement field~$j$.
	From equation~\eqref{eqn:eps_weak_fe}, if the displacement field is known,
	$\Bar{\eps}^{h}$ for an element is equal to:
	\begin{equation}
	\begin{split}
		q^{h}\Bar{\eps}^{h}
			&= q^{h}\eps^{h}_{\rm eq} 
				- \left. \frac{c^{2}}{h^{2}} \disc{q^{h}}\disc{\eps^{h}_{\rm eq}} 
					\right|_{j-1} 
				- \left. \frac{c^{2}}{h^{2}} \disc{q^{h}}\disc{\eps^{h}_{\rm eq}}
				\right|_{j} \\[1ex]
			&= q^{h} \brac{\eps^{h}_{\rm eq}  
			+ \left. \frac{c^{2}}{h^{2}} \disc{\eps^{h}_{\rm eq}} \right|_{j-1}
			- \left. \frac{c^{2}}{h^{2}} \disc{\eps^{h}_{\rm eq}} \right|_{j}}
	\end{split}
	\end{equation}	
	This is equivalent to:
	\begin{equation}
		\Bar{\eps}^{h} = \eps^{h}_{\rm eq} 
			+  \frac{c^{2}}{h} \brac{\left. u^{\prime \prime} \right|_{j}
				- 	\left. u^{\prime \prime} \right|_{j-1}}
	\end{equation}	
	which is a finite difference approximation of equation~\eqref{eqn:grad_strain},
	showing that the proposed variational formulation is identical to a
	finite-difference procedure in one-dimension for the case of 
	equally spaced nodal points.

\section{Numerical examples}
\label{sec:examples}
	Numerical examples is this section are intended to demonstrate 
	the objectivity of the formulation with respect to mesh refinement for
	strain softening problems, and to compare the computed results against a 
	known benchmark. It is well-known that classical, rate-independent continuum
	models are ill-posed when strain softening is introduced, which becomes
	evident in a severe sensitivity of the computed result to the spatial
	discretisation. One motivation for strain gradient dependent
	model is to provide regularisation in the presence of strain softening in
	order to avoid pathological mesh dependency.
	
	For all examples, the evolution of damage is given by
	equation~\eqref{eqn:damage_evolution}.
	The materials parameters are taken as: Young's modulus 
	$E = 20 \times 10^{3}$~MPa, $\kappa_{0}=0.0001$, $\kappa_{c} = 0.0125$ and
	$c=1$~mm. A Newton-Raphson procedure under displacement control is used to 
	solve the problem and the governing equations have been linearised 
	consistently.

\subsection{Objectivity with respect to spatial discretisation }
\begin{figure}
		\center\includegraphics[width=0.8\textwidth]{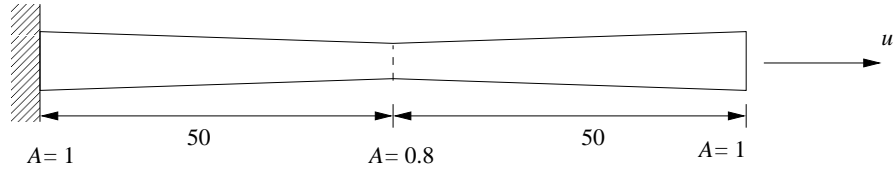}
	\caption{Linearly tapering bar.}
	\label{fig:tapered_bar}
	\end{figure}

\begin{figure}
		\center\includegraphics[scale=0.7]{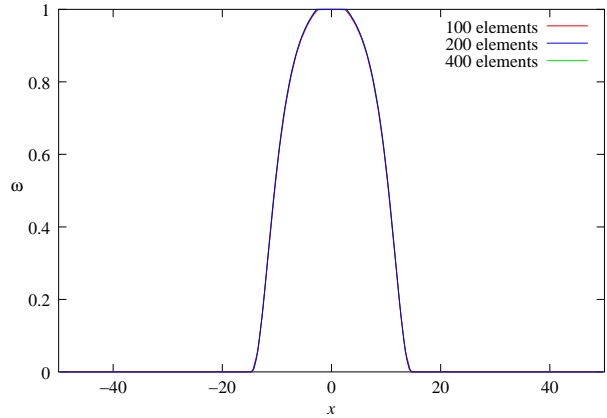}
	\caption{Damage profiles for the tapered bar.}
	\label{fig:tapered_dam}	
	\end{figure}

	The first test is for objectivity of the load--displacement response with
	respect to mesh refinement. A tapered bar (figure~\ref{fig:tapered_bar}) 
	is tested in tension. The bar has a cross-sectional area of one square unit
	at each end, and tapers linearly towards the centre where the area is
	0.8~square units. A displacement is applied incrementally at the right-hand 
	end.
	The response is examined for meshes with 100, 200 and 400 elements. 
	For each mesh, all elements are of equal size.
	\begin{figure}
		\center\includegraphics[scale=0.7]{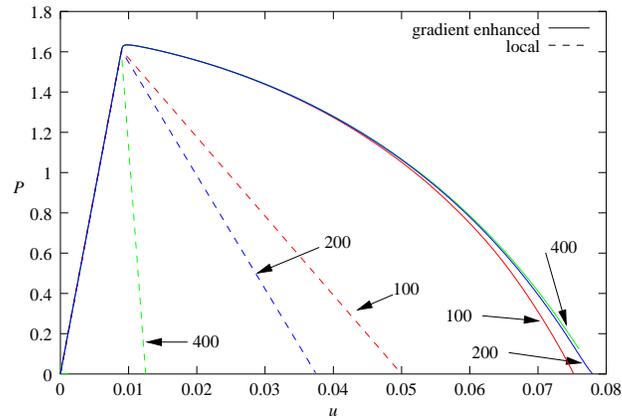}
	\caption{Load-displacement response for the tapered bar.}
	\label{fig:tapered_pd}	
	\end{figure}	
	Responses for the three meshes are shown in figure~\ref{fig:tapered_pd} 
	for both $c=1$ and $c=0$. Clearly, the introduction of strain gradient 
	effects has regularised 
	the problem, with the response for the three cases with $c=1$ being
	near identical. The response is further examined by comparing the damage
	profiles along the bar for the three regularised cases. The damage profiles,
	shown in figure~\ref{fig:tapered_dam}, are indistinguishable for the three
	meshes. 

\subsection{Comparison with a high-order of continuity numerical method}
\begin{figure}
		\center\includegraphics[width=0.8\textwidth]{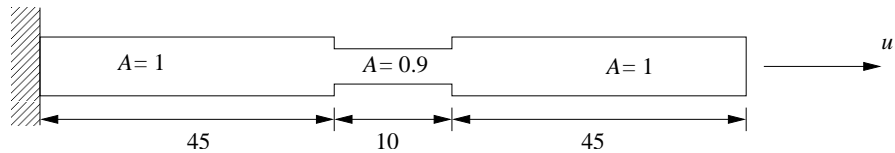}
	\caption{Bar with narrow central section.}
	\label{fig:bar_narrow}
	\end{figure}

	The second test involves a bar with a narrow section at the centre, 
	as shown in figure~\ref{fig:bar_narrow}. This problem was
	previously computed for the same strain gradient dependent damage model 
	using an element-free Galerkin method, which
	provides a high degree of continuity~\citep{askes:2000b}. 

\begin{figure}
		\center\includegraphics[scale=0.7]{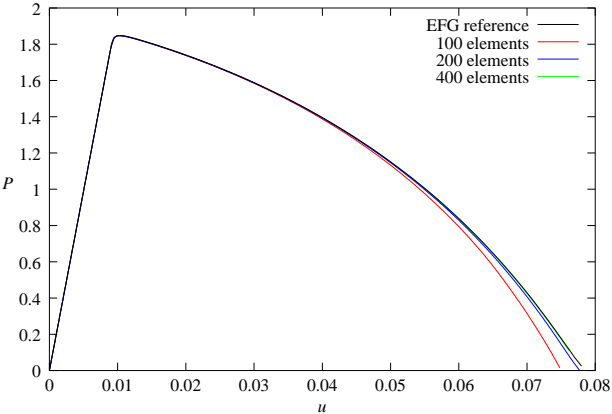}
	\caption{Load-displacement response.}
	\label{fig:bar_narrow_pd}	
	\end{figure}

\begin{figure}
		\center\includegraphics[scale=0.7]{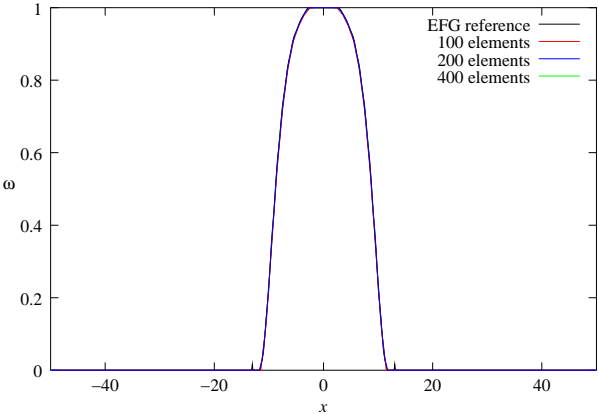}
	\caption{Damage profiles along the bar.}
	\label{fig:bar_narrow_dam}	
	\end{figure}

	This problem is computed using meshes with the same number of elements 
	as the previous example. For
	comparison, the computed load-displacement response from an element-free
	Galerkin method is also included~\citep{askes:2000b} for this problem. 
	It is clear from
	figure~\ref{fig:bar_narrow_pd} that the three meshes yield near-identical
	results and match the element-free Galerkin solution well.
	The damage profiles along the bar are shown in
	figure~\ref{fig:bar_narrow_dam}. The damage profile from 
	\citet{askes:2000b} is included as a reference.
	The computed results for all meshes are in excellent agreement with the
	benchmark.

\section{Conclusions}
	A discontinuous Galerkin  formulation has been developed for a strain
	gradient-dependent continuum model. The problem is split into two fields --
	the displacement and a deformation measure -- for generality.  The scalar
	field, which is a measure of the deformation, is dependent on gradients of 
	the strain field. Conventionally, this would require a $C^{1}$ finite
	element interpolation of the displacement field. By including element 
	interface terms in the Galerkin formulation, the need for high-order
	continuity is circumvented.
	
	The proposed formulation was tested for the simplest element configuration
	in one dimension -- piecewise continuous linear displacement and 
	discontinuous piecewise constant for the extra scalar field. For simple
	tests, the regularising properties of the strain gradient dependent model 
	were demonstrated
	and the results compared excellently with a benchmark result computed 
	using a numerical method with a high degree of continuity. These 
	preliminary results are
	promising and should be extended for higher-oder elements and 
	to multiple spatial dimensions. For the simple formulation adopted here,
	several terms in the weak from could be discarded. The importance of these
	terms must be assessed for higher-order interpolations.
	This work provides a first step towards a
	simple and well-founded finite element framework for modern strain gradient
	continuum models.

\section*{Acknowledgements}
	G.N.~Wells was partially supported by the J.~Tinsley Oden Faculty Research 
	Program,
	Institute for Computational Engineering and Sciences, The University of
	Texas at Austin for this work.
	The work of K.~Garikipati at University of Michigan was supported under NSF 
	grant CMS\#0087019.
	L.~Molari was supported under Progetto Marco Polo,  
	Universit\`{a} di Bologna. 
	The authors are grateful to R.L.~Taylor (University of California, Berkeley)
	for providing his program FEAP with discontinuous Galerkin capabilities.

\appendix

\section{Linearisation}
\label{append:linearisation}
	Effective solution of problems requires the consistent linearisation of the
	Galerkin problem. For the formulation, the fundamental unknowns are the 
	displacement
	$\vect{u}^{h}$ and~$\Bar{\eps}^{h}$. Linearisation requires 
	expressing the problem in terms of increments of the two unknowns.
	Taking the directional derivative of equation~\eqref{eqn:equil_cont_b},
	\begin{equation}
		\int_{\Omega} \gradn \vect{w}^{h} \ip \brac{1 - \omega}\mathcal{C} 
					\ip \Delta \strain \dif \Omega
		-\int_{\Omega} \gradn \vect{w}^{h} \ip \frac{\pd \omega}{\pd \Bar{\eps}}
			 \mathcal{C} 
					\ip \strain \Delta \Bar{\eps}^{h} \dif \Omega
		= \int_{\Gamma_{h}} \vect{w}^{h} \cdot \Delta \vect{h} \dif \Gamma 	
	\end{equation}
	where $\Delta \brac{\cdot}$ indicates a change in~$\brac{\cdot}$. 
	For brevity $\Delta \brac{\grad \vect{u}^{h}}$ is expressed 
	as~$\Delta \strain^{h}$. Since the gradient is a linear operator, $\Delta
	\brac{\grad \vect{u}} = \grad\brac{\Delta \vect{u}}$.
	Similarly, equation~\eqref{eqn:eps_weak_b} is linearised by
	taking the directional derivative,
	\begin{multline}
		\int_{\Omega} q^{h} \Delta \Bar{\eps}^{h} \dif \Omega
	   -\int_{\Omega} q^{h} \frac{\pd \eps_{\rm eq}}{\pd \strain}
	    			\ip \Delta \strain^{h} \dif \Omega
	   +\int_{\Tilde{\Omega}} \gradn q^{h} \cdot c^{2} 	    	  
			 		 \gradn \brac{\frac{\pd \eps_{\rm eq}}{\pd \strain}
					\ip \Delta \strain^{h}} \dif \Omega \\
	   -\int_{\Tilde{\Gamma}} \disc{q^{h}} \cdot c^{2} 
	    	 \ave{ \gradn\brac{\frac{\pd \eps_{\rm eq}}{\pd \strain} \ip 
			 \Delta \strain^{h}}} \dif \Gamma
	   -\int_{\Tilde{\Gamma}} \ave{\gradn q^{h}} \cdot c^{2} 
	    	 \disc{ \frac{\pd \eps_{\rm eq}}{\pd \strain} \ip \Delta
			\strain^{h}} \dif \Gamma  \\
	   +\int_{\Tilde{\Gamma}} \frac{c^{2}}{h_{e}} \disc{q^{h}} \cdot 
	    	 \disc{ \frac{\pd \eps_{\rm eq}}{\pd \strain} \ip \Delta
			\strain^{h}} \dif \Gamma
		  = \int_{\Gamma} q^{h} \Delta \epsbc \dif \Gamma.  	
	\end{multline}

\section{Finite element formulation}
\label{append:fe}
	The finite element formulation  is elaborated here for the case of 
	piecewise continuous linear $\vect{u}^{h}$ and piecewise constant 
	$\Bar{\eps}^{h}$.
	It can be extended to the more general case
	of arbitrary interpolation orders.
	
	Formulation of the stiffness matrix consists of two keys steps. 
	The first is the usual loop over all elements. This yields a stiffness 
	matrix for each element $\vect{k}_{e}$ of the form
	\begin{equation}
		\vect{k}_{e} = 
			\begin{bmatrix}
				\vect{k}_{uu} 	& \vect{k}_{u\Bar{\eps}}  \\
				\vect{k}_{\Bar{\eps}u} & \vect{k}_{\Bar{\eps}\Bar{\eps}}
			\end{bmatrix}
	\end{equation}
	where the components of the matrix $\vect{k}_{e}$ are 
	\begin{align}
		\vect{k}_{uu} 
			&= \int_{\Omega_{e}} \brac{1-\omega}\trans{\vect{B}} \vect{D} \vect{B}
					\dif \Omega				\\ 
		\vect{k}_{u\Bar{\eps}} 
			&=	- \int_{\Omega_{e}} \trans{\vect{B}} \frac{\pd \omega}{\pd
			\Bar{\eps}} \vect{D} \strain \vect{N}_{\Bar{\eps}} 
						\dif \Omega			\\ 
		\vect{k}_{\Bar{\eps}u} 
			&= -\int_{\Omega_{e}} \trans{\vect{N}}_{\Bar{\eps}}  
				\trans{\brac{\frac{\pd \eps_{\rm eq}}{\pd \strain}}}
				 \vect{B} \dif \Omega \\ 
		\vect{k}_{\Bar{\eps}\Bar{\eps}} 
			&= \int_{\Omega_{e}} \trans{\vect{N}}_{\Bar{\eps}} 
					\vect{N}_{\Bar{\eps}} \dif \Omega				
	\end{align}
	where $\vect{B}$ is the usual finite element matrix containing spatial
	derivatives of the shape functions related to the displacement field,
	$\vect{D}$ is the elastic constitutive tensor in matrix form and
	$\vect{N}_{\Bar{\eps}}$ contains the shape functions relating the
	$\Bar{\eps}$. The strain is expressed in engineering column vector format.
	Once formed, an element element stiffness matrix is 
	assembled into the global system of equations as usual.

	The next, non-standard, step is a loop over all element interfaces. For
	this, `information' is required for both the elements that are connect to
	the interface. The stiffness matrix at the interface two equal-order
	elements is twice the size of the stiffness matrix of a single element. It
	can be expressed as: 
	\begin{equation}
		\vect{k}_{i} = 
			\begin{bmatrix}
				\vect{k}_{u_{1}u_{1}} 	& \vect{k}_{u_{1}\Bar{\eps}_{1}}
				  &  \vect{k}_{u_{1}u_{2}} 	& \vect{k}_{u_{1}\Bar{\eps}_{2}} \\
				\vect{k}_{\Bar{\eps}_{1}u_{1}} 
						& \vect{k}_{\Bar{\eps}_{1}\Bar{\eps}_{1}} &
				\vect{k}_{\Bar{\eps}_{1}u_{2}} 
						& \vect{k}_{\Bar{\eps}_{1}\Bar{\eps}_{2}} \\ 
				\vect{k}_{u_{2}u_{1}} 	& \vect{k}_{u_{2}\Bar{\eps}_{1}}
				  &  \vect{k}_{u_{2}u_{2}} 	& \vect{k}_{u_{2}\Bar{\eps}_{2}} \\
				\vect{k}_{\Bar{\eps}_{2}u_{1}} 
						& \vect{k}_{\Bar{\eps}_{2}\Bar{\eps}_{1}} &
				\vect{k}_{\Bar{\eps}_{2}u_{2}} 
						& \vect{k}_{\Bar{\eps}_{2}\Bar{\eps}_{2}}  
			\end{bmatrix}
	\end{equation}
	where the subscripts `1' and `2' denote the element on either side of the
	surface. For the case of linear $\vect{u}^{h}$ and constant $\Bar{\eps}$, 
	only the terms $\vect{k}_{\Bar{\eps}_{j}u_{k}}$ are non-zero. It is equal
	to:
	\begin{equation}
		\vect{k}_{\Bar{\eps}_{j}u_{k}} 
			= \int_{\Tilde{\Gamma}_{i}} \frac{c^{2}}{h_{e}}
				\trans{\vect{N}}_{\Bar{\eps}j} \trans{\vect{n}}_{j}
				\vect{n}_{k}
				\trans{\brac{\frac{\pd \eps_{\rm eq}}{\pd \strain}}}		
				\vect{B}_{k}		
						\dif \Gamma
	\end{equation}
	where the indices $j$ and $k$ run from one to two, corresponding to sides of
	the interface. 	Note that in the usual case of 
	$\vect{n}_{1} = - \vect{n}_{2}$,
	$\trans{\vect{n}}_{i} \vect{n}_{j} = 1$ if $i=j$, and  
	$\trans{\vect{n}}_{i} \vect{n}_{j} = -1$ if~$i\ne j$.

\bibliographystyle{cmame}
\bibliography{papers_full,papers,papers_gradient,dg}
%
%
\end{document}